\documentclass{amsproc}

\newtheorem{theorem}{Theorem}[section]

\newtheorem{corollary}[theorem]{Corollary}

\theoremstyle{definition}
\newtheorem{definition}[theorem]{Definition}
\newtheorem{example}[theorem]{Example}

\theoremstyle{remark}
\newtheorem{remark}[theorem]{Remark}
\newtheorem*{acknowledgments}{Acknowledgments}

\numberwithin{equation}{section}

\newcommand{\nn}{{\newline}}
\begin{document}

\title{An analogue of Bratteli-Jorgensen loop group actions for GMRA's}

\author{L. W. Baggett}
\address{Department of Mathematics, Campus Box 395, University of Colorado, Boulder, CO, 80309-0395}
\email{baggett@euclid.colorado.edu}
\author{P. E. T. Jorgensen}
\address{Department of Mathematics, University of Iowa, 14 MacLean Hall, 
Iowa City, IA, 52242-1419}
\email{jorgen@math.uiowa.edu}
\thanks{The first two named authors were supported by a US-NSF 
Focused Research Group (FRG) grant.}
\author{K. D. Merrill}
\address{Department of Mathematics, Colorado College, Colorado Springs, CO 80903-3294}
\email{kmerrill@coloradocollege.edu}
\author{J. A. Packer}
\address{Department of Mathematics, Campus Box 395, University of Colorado, Boulder, CO, 80309-0395}
\email{packer@euclid.colorado.edu}
\date{June 9, 2003}


\begin{abstract}
Several years ago, O. Bratelli and P. Jorgensen developed the concept of $m$-systems of filters for dilation by a positive integer $N>1$ on $L^2(\mathbb R)$.  They constructed a loop group action on $m$-systems.  By work of Mallat and Meyer, these $m$-systems are important in constructing  multi-resolution analyses and wavelets associated to dilation by $N$ and translation by $\mathbb Z$ on $L^2(\mathbb R)$. In this paper, we discuss an extension of this loop-group construction to generalized filter systems, which we will call ``$M$-systems," associated with generalized multiresolution analyses.  In particular, we show that every multiplicity function has an associated generalized loop group which acts freely and transitively on the set of $M$-systems corresponding to the multiplicity function.  The results of Bratteli and Jorgensen correspond to the case where the multiplicity function is identically equal to $1.$
\end{abstract}

\maketitle


\section{Introduction}
In work that first appeared in the late 1990's and has since been elaborated upon in their book \cite{BJ}, O. Bratteli and P. Jorgensen related filter functions corresponding to a multiresolution analysis  for dilation by an integer $N>1$ to representations of the Cuntz algebras.  Brattelli and Jorgensen obtained results about frame wavelet families and orthonormal wavelet families arising from multiresolution analyses using techniques from operator theory, and in addition they were able to construct a free and transitive action of a ``loop group" on ``$\text{Lip}_1\;m$- systems," consisting of families of H\"{o}lder-continuous filter functions for dilation by $N$.  

At the same time, L. Baggett, H. Medina and K. Merrill developed a theory of generalized multiresolution analyses for dilation by $N$ \cite{BMM}, which we abbreviate here using the acronym ``GMRA".  Somewhat later, in \cite{B}, L. Baggett, J. Courter and K. Merrill associated to these objects generalized filter functions, which in this paper we will generalize to the concept of ``$M$- systems".
Under appropriate conditions, Baggett, Courter and Merrill were able to construct a GMRA and an associated frame wavelet system by using the $M$-systems.  In both of these papers, the so-called ``multiplicity function" $\mu$ associated to a GMRA and the ``consistency equation" that the multiplicity function satisfies turned out to be of great importance in the proofs of theorems, and that will be the case here as well.   

In this paper, we will review the loop group action for ``classical $m$- systems." We will then define the notion of $M$-systems corresponding to a fixed multiplicity function $\mu$. The $M$-systems turn out to be a set of Borel cross-sections of a Borel vector bundle over a certain subset of 
$\bigsqcup_{i=1}^c\mathbb T,$ where $c$ the essential supremum of $\mu$.  
The generalized filter systems of Baggett, Courter and Merrill (see \cite{B})turn out to be prime examples of $M$-systems, and of course motivate the general definition. Our main result is the construction of a generalized loop group corresponding to a multiplicity function $\mu.$  We describe this group here as the group of sections of a group bundle over $\mathbb T.$ 
  We will show that this generalized loop group acts freely and transitively on the set of $M$-systems.  
 
For the purposes of this article, we will for the most part regard a family of $M$-systems corresponding to a given multiplicity function as an abstract set which is acted on by the loop group.  However, our ultimate aim, which we plan to present in a later work, is to use the operator approach pioneered by Bratteli and Jorgensen in the classical case to weaken the conditions given in \cite{B} under which an $M$-system gives rise to a GMRA, and hence to a normalized tight frame wavelet family in $L^2(\mathbb R)$ for dilation by $N.$  This would provide another approach to the construction of normalized tight frame wavelets. 

\section{Preliminaries: the classical case, and a loop group action on $m$-systems}

We first review the notion of a normalized tight frame in a Hilbert space.
\begin{definition} 
  A sequence $\{ x_n : n \in {\mathbb N} \}$ of elements in a Hilbert
  space $\mathcal H$ is said to be a {\it frame} if there are real
  constants $C,D > 0$ such that
  \begin{equation} \label{ineq-frame}
      C \cdot \|x,x\|^2\leq
      \sum_{n=1}^{\infty} |\langle x,x_j \rangle|^2  \leq
      D \cdot \| x,x \|^2
  \end{equation}
  for every $x \in \mathcal H$.  If $C=D=1,$ we say the set $\{ x_n : n \in {\mathbb N} \}$ is a {\it normalized tight frame} (abbreviated NTF) for ${\mathcal H}.$
\end{definition}

It is possible to show that a subset $\{ x_n : n \in {\mathbb N} \}$ of a Hilbert space ${\mathcal H}$ is a NTF for ${\mathcal H}$ if and only if for every $x\;\in\;{\mathcal H}$ the following reconstruction formula is satisfied:
\begin{equation} \label{reconform}
      x\;=\; \sum_{n=1}^{\infty} \langle x,x_n \rangle x_n.
  \end{equation}

Any orthonormal basis for ${\mathcal H}$ is a normalized tight frame, but not conversely, as frames can have redundancy which does not occur in the orthonormal basis case.

We now move on to the case where ${\mathcal H}\;=\;L^2(\mathbb R)$ and review the definition of a normalized tight frame wavelet family in this context. Fix an integer $N>1.$ 
Fix an integer $N>1.$ In the standard way, we define dilation and translation operators on $L^2(\mathbb R)$ by 
$$D(f)(t)\;=\;\sqrt{N}f(Nt),$$
$$T(f)(t)\;=\;f(t-1),\;f\in\;L^2(\mathbb R).$$

As this paper is mainly concerned with low and high pass filters,  it is natural to define orthonormal wavelets and tight frame wavelets in the frequency domain.
Let ${\mathcal F}$ denote the Fourier transform on $L^2(\mathbb R),$ defined by 
$${\mathcal F}(f)(x)\;=\;\int_{\mathbb R}f(t)e^{-2\pi ixt}dt.$$  
Set 
$${\widehat{D}}={\mathcal F}D{\mathcal F}^*,\;\widehat{T}={\mathcal F}T{\mathcal F}^*.$$
Then 
$${\widehat{D}}(f)(x)=(\sqrt{N})^{-1}f(N^{-1}x),$$ and 
$$\widehat{T}(f)(x)=e^{-2\pi ix}f(x),\;f\in L^2(\mathbb R).$$
\begin{definition}
A {\bf {normalized tight frame wavelet family in the frequency domain}} for dilation by $N>1$ is a subset $\{\Psi_1,\cdots,\Psi_m\}\subseteq L^2(\mathbb R)$ such that
$$\{{\widehat{D}}^j \widehat{T}^v(\Psi_i):\;1\leq i \leq m,\;j\in\mathbb Z,\;v\in\mathbb Z\}$$
is an NTF for $L^2(\mathbb R).$ If the $\{{\widehat{D}}^j {\widehat {T}}^v(\Psi_i):\;1\leq i \leq m,\;j\in\mathbb Z,\;v\in\mathbb Z\}$ form an orthonormal basis for $L^2(\mathbb R),$ the family $\{\Psi_1,\cdots,\Psi_m\}$ is called an orthonormal wavelet family for dilation by $N.$
\end{definition}

A fundamental algorithm for constructing wavelet families is the concept of multiresolution analysis (MRA) developed by S. Mallat and Y. Meyer \cite{Ma}, and key tools for constructing the MRA's are filter functions for dilation by $N.$ 
\begin{definition} 
Let $N$ be a positive integer greater than $1.$  A {\bf low-pass filter} $m_0$ for dilation by $N$ is a $\mathbb Z$-periodic function 
\nn $m_0:\;\mathbb R\;\rightarrow\;\mathbb C$ which satisfies the following conditions:
\renewcommand{\labelenumi}{(\roman{enumi})}
\begin{enumerate}
\item $m_0(0)\;=\;\sqrt{N}$ (``low-pass condition")
\item $\sum_{l=0}^{N-1}|m_0(x+l/N)|^2=N$ a.e.;
\item $m_0$ is H\"{o}lder continuous at $0;$ 
\item $m_0$ is  non-zero in a sufficiently large neighborhood of $0$ (``Cohen's condition", c.f. \cite{Ch}).
\end{enumerate}
\end{definition}
Sometimes in the above definition, condition (iv), Cohen's condition, is dropped.

Given a low-pass filter $m_0$ for dilation by $N,$ then there is a canonical way to construct a ``scaling function", or ``father wavelet" associated to the filter. We set
$$\Phi(x)\;=\Pi_{i=1}^{\infty}[\frac{m_0(N^{-i}x)}{\sqrt{N}}].$$
Then $\Phi$ converges a.e. and is an element of $L^2(\mathbb R).$ We call $\Phi$ a {\bf  scaling function in the frequency domain} for dilation by $N.$  

To use the low-pass filter and the associated scaling function to construct a wavelet family for dilation by $N$, we need $N-1$ other functions that are high pass filters in the sense of following  definition.  
\begin{definition}
Let $N$ be a positive integer greater than $1,$ and let  $m_0$  be a  low-pass filter for dilation by $N.$ A set of essentially bounded  measurable $\mathbb Z$ periodic functions $m_1,m_2,\cdots,m_{N-1}$ defined on $\mathbb R$ are called {\bf high-pass filters} associated to $m_0,$ 
if 
$$\sum_{l=0}^{N-1}{\overline{m_i(x+\frac{l}{N})}}m_j(x+\frac{l}{N})\;=\;\delta_{i,j}N,\;0\leq\;i,j\;\leq N-1.$$
\end{definition}
We can express the filter conditions more concisely as follows: 
the 
$N\times N$ matrix-valued function on $x\in\mathbb R/\mathbb Z$ 
\begin{equation}
\label{matrixfilt}
x\;\mapsto\;(\frac{m_j(x+\frac{l}{N})}{\sqrt{N}})_{0\leq j,l \leq N-1}
\end{equation}
is a unitary matrix-valued function $\forall x\;\in\;\mathbb T.$
\begin{definition}\cite{BJ}
Let $m_0$ be a $\text{Lip}_1$ (i.e. H\"{o}lder continuous) low-pass filter for dilation by $N,$ which need not satisfy Cohen's condition. Suppose that 
 $m_1,m_2,\cdots,m_{N-1}$ are a collection of $N-1\;\text{Lip}_1$ high-pass filters associated to $m_0.$ The $N$-tuple $(m_0,m_1,\cdots, m_{N-1})$ is called a {\bf classical $\text{Lip}_1\;m$-system} for dilation by $N$.
\end{definition}
\begin{theorem} (\cite{Ma}, \cite{Me}, \cite{L}, \cite{BJ})  
Let $N$ be a positive integer greater than $1,$ let 
$(m_0,m_1,\cdots, m_{N-1})$ be a classical $m$-system for dilation by $N,$
 and let 
$\Phi$ be the scaling function constructed from $m_0$ as above.Then 
$$\{ \Psi_1=\widehat{D}(m_1\Phi),\;\Psi_2=\widehat{D}(m_2\Phi),\;\cdots,\;\Psi_{N-1}=\widehat{D}(m_{N-1}\Phi)\}$$ is an NTF wavelet family in the frequency domain for dilation by $N.$ If Cohen's condition is satisfied, the $\{\Psi_k\}$'s  form an orthonormal wavelet family. 

\end{theorem}
\begin{remark} 
In their proof of this result, Bratteli and Jorgensen used a representation of the Cuntz algebra ${\mathcal O}_N$ arising from the $m$-system.
\end{remark}
So we see that classical $m$-systems are very useful in the construction of NTF wavelet families. 

In their book \cite{BJ}, Bratteli and Jorgensen constructed a so-called ``loop group" which acts freely and transitively on the set of $m$-systems.  We modify the conventions slightly in the definition below, as it will be more easy to adapt to the GMRA $M$-system case. 
\begin{definition}
Fix a positive integer $N>1.$ We say a matrix valued function $K:\mathbb R/\mathbb Z:\;\rightarrow\;U(N,\mathbb C)$ is in $\text{Lip}_1(MF,N)$ if 
\renewcommand{\labelenumi}{(\roman{enumi})}
\begin{enumerate}
\item Each component function $k_{i,j}$ is $Lip_1,\;0\leq i,j \leq\;N-1,$
\item $K(0)=Id_{N\times N}.$
\end{enumerate}
\end{definition}
Then the main result of Bratteli and Jorgensen about the loop group is as follows:
\begin{theorem}\cite{BJ} 
The loop group $\text{Lip}_1(MF,N)$ acts freely and transitively on $\text{Lip}_1\;m$ systems for dilation by $N$ via the formula 
$$K\cdot(m_0,m_1,\cdots,m_{N-1})^T](x)\;=\;K(Nx)(m_0(x),m_1(x),\cdots,m_{N-1}(x))^T.$$
\end{theorem}
\begin{proof}  
We sketch the proof of this result for completeness. We first prove transitivity.  If $(m_0,m_1,\cdots,m_{N-1})$ and $(\widetilde{m}_0,\widetilde{m}_1,\cdots,\widetilde{m}_{N-1})$ are two $\text{Lip}_1\;m$ systems for dilation by $N,$ then defining $K(x)\;=\;(k_{i,j}(x))$ by 
$$k_{i,j}(x)\;=\;\frac{1}{N}\sum_{l=0}^{N-1}\widetilde{m}_i(\frac{x+l}{N})\overline{m_j(\frac{x+l}{N})},$$ 
one calculates that 
\begin {eqnarray}
[K(Nx)(m_0,(x)m_1(x),\cdots,m_{N-1}(x))^T]_i &=&\sum_{j=0}^{N-1}k_{i,j}(Nx)m_j(x) \nonumber \\
&=&\sum_{j=0}^{N-1}\frac{1}{N}\sum_{l=0}^{N-1}\widetilde{m}_i(\frac{Nx+l}{N})\overline{m_j(\frac{Nx+l}{N})}m_j(x) \nonumber \\
&=&\sum_{l=0}^{N-1}\widetilde{m}_i(x+\frac{l}{N})[\frac{1}{N}\sum_{j=0}^{N-1}\overline{m_j(x+\frac{l}{N})}m_j(x)] \nonumber \\
&=&\sum_{l=0}^{N-1}\widetilde{m}_i(x+\frac{l}{N})\delta_{l,0}\;=\;\widetilde{m}_i(x), \nonumber
\end{eqnarray}
by the conditions on the filter systems coming from the fact that the matrix described in Equation \ref{matrixfilt} is unitary.
Thus $$K\cdot(m_0(x),m_1(x),\cdots,m_{N-1}(x))^T\;=\;[(\widetilde{m}_0,\widetilde{m}_1,\cdots,\widetilde{m}_{N-1})(x)]^T,$$ and
$\text{Lip}_1(MF,N)$ acts freely on the set of $m$-systems.

As for freeness, let $(m_0,m_1,\cdots,m_{N-1})$ be a $m$-system.  Recall that we can associate to this $m$-system the $N\times N$ matrix whose entries are in $\text{Lip}_1(\mathbb T)$ which is unitary $\forall x\in\mathbb T$ given by 
$${\mathcal M}(x)\;=\;(\frac{m_i(x+\frac{j}{N})}{\sqrt{N}})_{0\leq i,j \leq N-1}.$$
Suppose that there exists $K\in\text{Lip}_1(MF,N)$ satisfying   
$$K\cdot[m_0(x),m_1(x),\cdots,m_{N-1}(x)]^T\;=\;[m_0(x),m_1(x),\cdots,m_{N-1}(x)]^T,$$
for some $m$ system $(m_0,m_1,\cdots,m_{N-1}),$ for all $x\;\in\;\mathbb T.$
But then for all $x\;\in\;\mathbb T$ and for $0\;\leq\; j\;\leq\; N-1$ we have 
$$K \cdot[m_0(x+\frac{j}{N}),m_1(x+\frac{j}{N}),\cdots,m_{N-1}(x+\frac{j}{N})]^T\;=\;[m_0(x+\frac{j}{N}),m_1(x+\frac{j}{N}),\cdots,m_{N-1}(x+\frac{j}{N})]^T,$$
i.e. $$K(N(x+\frac{j}{N}))[m_0(x+\frac{j}{N}),m_1(x+\frac{j}{N}),\cdots,m_{N-1}(x+\frac{j}{N})]^T\;=$$
$$[m_0(x+\frac{j}{N}),m_1(x+\frac{j}{N}),\cdots,m_{N-1}(x+\frac{j}{N})]^T,$$
for $0\;\leq\; j\;\leq\; N-1.$
But, since the entries of $K$ are $\mathbb Z$-periodic, this means (in terms of matrix multiplication) exactly that  
$$K(Nx){\mathcal M}(x)\;=\;{\mathcal M}(x),\; \forall x\in\mathbb T,$$
where ${\mathcal M}(x)=(m_i(x+j/N)/\sqrt{N})_{0\leq i,k \leq N-1}$ is the unitary-valued matrix function defined in Equation \ref{matrixfilt}. 
We thus obtain $K(Nx)\;=\;{\mathcal M}(x){\mathcal M}(x)^{\ast}\; \forall x\in\mathbb T.$ Thus $K(x)$ is the identity matrix for every $x\in\mathbb T,$ as we desired to show.
\end{proof} 
\section{Generalized Multiresolution Analyses and Generalized Filter Systems:} 

Not all wavelet families come  come from MRA's and $m$-systems. There are various ways to approach the general case, c.f. \cite{BL}, \cite{HLPS}. In 1999, L. Baggett, H. Medina and K. Merrill developed the theory of generalized multiresolution analyses (GMRA's) to deal with this more general setting \cite{BMM}.  Again to simplify matters we limit our definition to the frequency domain for now.
\begin{definition} \cite{BMM} A {\bf generalized multiresolution analysis} (GMRA) for dilation by $N$ is a sequence $\{\widehat{V_i}\}_{i\in\mathbb Z}$ of closed subspaces of $L^2(\mathbb R)$ satisfying the following conditions:  
\renewcommand{\labelenumi}{(\roman{enumi})}
\begin{enumerate}
\item $\;\;\cdots \widehat{V_{-1}}\subseteq\widehat{V_0}\subseteq\widehat{V_1}\cdots$ (the $\widehat{V_i}$ are nested) 
\item $(\widehat{D})^i(\widehat{V_0})\;=\;\widehat{V_i},\;i\in\mathbb Z;$
\item $\overline{\cup_{i\in\mathbb Z}\widehat{V_i}}\;=\;L^2(\mathbb R),\;\cap_{i\in\mathbb Z}\widehat{V_i}=\{0\};$
\item $\widehat{V_0}$ is invariant under all powers of $\widehat{T}.$
\end{enumerate}
\end{definition}
Baggett, Merrill and Medina were able to develop several characteristic invariants associated to a GMRA using spectral theory.  We summarize their key results in the following theorem.
\begin{theorem} \cite{BMM} Given a GMRA in $L^2(\mathbb R)$ corresponding to dilation $\widehat{D}$ by $N$ and the transform of integer translation $\widehat{T},$ there is a unique sequence of Borel subsets $S_1\supseteq S_2\supseteq\cdots$ of $\mathbb T$ and a unitary operator $J:\widehat{V_0}\;\rightarrow\;\oplus_j\;L^2(S_j)$ such that 
$$[J(T(f)]_j(x)\;=\;e^{2\pi ix}[J(f)]_j(x).$$
The function $\mu(x)\;=\;\sum_j\chi_{S_j}(x)$ defined on $\mathbb T$ is called the {\bf multiplicity function} corresponding to the GMRA $\{\widehat{V_i}\}_{i\in\mathbb Z},$ and satisfies 
$$\mu(x)\;\leq\;\sum_{l=0}^{N-1}\mu(\frac{x+l}{N})\;\text{a.e.}.$$
\end{theorem}
If $\mu$ is essentially bounded, with $c\;=\text{ess sup}\;\mu(x),$ it is possible to define the conjugate multiplicity function 
$$\tilde{\mu}(x)\;=\;\sum_{l=0}^{N-1}\mu(\frac{x+l}{N})-\mu(x)$$
By definition, $\mu$ and $\tilde{\mu}$ satisfy the so-called ``consistency equation":
\begin{equation}
\label{coneq}
\mu(x)\;+\;\tilde{\mu}(x)\;=\;\sum_{l=0}^{N-1}\mu(\frac{x+l}{N}).
\end{equation}
@3 43nq4j 5gq5 5g3 dqw3 $\mu\equiv 1$ corresponds to the MRA case, and the case $\tilde{\mu}\equiv = 1$ corresponds to the case where there is a single orthonormal wavelet.
\newline Note also by the definition of $\mu$ that the sets $\{S_i\}_{i=1}^c$ satisfy the following condition:
\begin{equation}
\label{Seq}
S_i\;=\;\{x\in\;\mathbb T\; | \;\mu(x)\geq\;i\}.
\end{equation}
Sets $\widetilde{S}_i$ are defined analogously to the sets $S_i$ by $\widetilde{S}_i\;=\;\{x\in\;\mathbb T\;|\;\tilde{\mu}(x)\geq\; i\}.$  If $d\;=\text{ess sup}\;\tilde{\mu}(x),$ we have 
$$\widetilde{S}_1\;\supseteq\;\widetilde{S}_2\;\cdots\;\supseteq\;\widetilde{S}_d.$$ 

Baggett, Courter and Merrill in \cite{C2} and \cite{B} then generalized the Mallat and Meyer algorithm for constructing wavelets from filters to this GMRA setting.  They first generalized the concept of low and high pass filters.  Suppose $\mu$ is a Borel integer-valued function that is essentially bounded by c on $\mathbb T,$ is constant in a neighborhood of the origin, and satisfies technical conditions that allow it to be a multiplicity function for a GMRA (c.f. \cite{BM}). They then defined ``generalized conjugate mirror filters", which correspond to low pass filters in the classical case, to be functions $\{h_{i,j}\}_{1\leq\;i,j\;\leq c},$ where each $h_{i,j}$ is supported on $S_j,$ and the following orthogonality condition holds:
\begin{equation}
\label{ortho1}
\sum_{j=1}^c\sum_{l=0}^{N-1}h_{i,j}(\frac{x+l}{N})\overline {h_{k,j}(\frac{x+l}{N})}\;=\;N\delta_{i,k}\chi_{S_i}(x).
\end{equation}
Similarly, they defined ``complementary conjugate mirror filters," an analogue of classical high pass filters, to be functions
$\{ g_{k,j}\}_{1\leq\;k\;\leq d,\;1\leq\;j\;\leq c},$ where each $g_{k,j}$ is supported on $S_j$, and 
\begin{equation}
\label{ortho2}
\sum_{j=1}^c\sum_{l=0}^{N-1}g_{k,j}(\frac{x+l}{N})\overline {g_{k',j}(\frac{x+l}{N})}\;=\;N\delta_{k,k'}\chi_{\widetilde{S}_k}(x),
\end{equation}
and
\begin{equation}
\label{ortho3}
\sum_{j=1}^c\sum_{l=0}^{N-1}h_{i,j}(\frac{x+l}{N})\overline {g_{k,j}(\frac{x+l}{N})}=0,\;\forall\;i,\;k.
\end{equation}
Examples of functions $g_{k,j}$ and $h_{i,j}$ that satisfy these conditions can be built by an explicit algorithm, and then modified to mimic examples of classical filters \cite{B}.  

Under appropriate conditions on these generalized filter functions, Baggett, Courter and Merrill then used them to construct a finite tight frame wavelet family $\{\Psi_1,\cdots,\Psi_d\}\subseteq L^2(\mathbb R)$. 
They first built generalized scaling functions 
\newline
$\{\Phi_1,\Phi_2, \cdots, \Phi_c\}\;\subseteq\;\widehat{V_0}$ using an infinite product construction involving dilates of a matrix with periodizations of $\{h_{i_j}\}$ as entries (c.f. Theorem 3.4 \cite{B}). The $\{\Phi_i\}_{i=1}^c$ appear as the first column in the infinite product matrix. 

Given the above notation and construction, they then have:
\begin{theorem} \cite{B} Let  $\{(h_{i,j}\}_{1\leq\;i,j\;\leq c}$ and $\{ g_{k,j}\}_{1\leq\;k\;\leq d,\;1\leq\;j\;\leq c}$ be generalized filter functions associated to the multiplicity function
$\mu,$ that satisfy appropriate conditions and let $\{\Phi_1,\Phi_2, \cdots, \Phi_c\}\;\subseteq\;\widehat{V_0}$ be generalized scaling functions constructed as described above.  Setting 
$$\Psi_k\;=\widehat{D}(\sum_{j=1}^c\;g_{k,j}\Phi_j),\;1\leq\;k\;\leq d,$$
the $\{\Psi_k\;|\;1\leq\;k\;\leq d\}$ form a NTF wavelet family for dilation by $N.$
\end{theorem}
As mentioned above, our ultimate aim is to use the approach of Bratteli and Jorgensen to weaken the restrictive technical conditions required in this theorem (which we do not make explicit here) to hypotheses more in line with the classical result.  To move toward that goal, we now use these generalized filter functions to define precise analogues of the classical high and low pass filters. In the next section, these definitions will allow us to define generalized M systems and develop the loop group action on them.   
\begin{definition} 
Given a multiplicity function $\mu$ that is constant in a neighborhood of the origin, with its associated sequence of sets, $\{S_i\:|\;1\leq\; i\;\leq c\}$, a {\bf generalized low-pass filter}  for dilation by $N$ is a collection of functions $\{h_{i,j}\}_{1\leq\;i,j\;\leq c},$ that satisfy the following conditions:
\renewcommand{\labelenumi}{(\roman{enumi})}
\begin{enumerate}
\item $h_{i,j}$ is supported on $S_j$;
\item $h_{i,j}(0)\;=\;\sqrt{N}\delta_{1,j}$ (``low-pass condition");
\item $\{h_{i,j}\}$ satisfy the orthogonality condition \ref{ortho1};
\item $h_{i,j}$ is $\text{Lip}_1$ at $0.$ 
\end{enumerate}
\end{definition}
\begin{definition} 
Let $\{h_{i,j}\}_{1\leq\;i,j\;\leq c}$ be a generalized low-pass filter associated with the multiplicity function $\mu$.  Suppose that both $\mu$ and the conjugate multiplicity function $\tilde{\mu}$ are constant in a neighborhood of the origin.   Define the sets $\{\widetilde{S}_k\:|\;1\leq\; k\;\leq d\}$ as above.  An associated {\bf generalized high-pass filter} is a collection of functions $\{ g_{k,j}\}_{1\leq\;k\;\leq d,\;1\leq\;j\;\leq c}$ that satisfy the following conditions:
\renewcommand{\labelenumi}{(\roman{enumi})}
\begin{enumerate}
\item $g_{k,j}$ is supported on $S_j,$;
\item $\{g_{k,j}\}$ satisfy the orthogonality conditions \ref{ortho2} and \ref{ortho3};
\item $g_{k,j}$ is $\text{Lip}_1$ at $0.$ 
\end{enumerate}
\end{definition}

\section{A Generalized Loop Group Action on the Generalized Filter Systems}

Let $\{h_{i,j}\}_{1\leq\;i,j\;\leq c}$ and 
$\{ g_{k,j}\}_{1\leq\;k\;\leq d,\;1\leq\;j\;\leq c}$ be generalized low-pass and high-pass filter functions defined as in the previous section.  
Since $\oplus_{i=1}^c\;L^2(S_i)\;\cong\;L^2(\bigsqcup_{i=1}^c S_i),$ we can suppress the second index of the filter functions and view generalized filter functions as a vector ($c+d$-tuple) of functions: 
$$(h_1,h_2,\;\cdots\;h_c,g_1,g_2,\cdots,g_d)\;\in\;[L^2(\bigsqcup_{i=1}^c S_i)]^{c+d}.$$
Further, we note that for any fixed $x$, the output of the vector of functions 
\newline $(h_1,h_2,\;\cdots\;h_c,g_1,g_2,\cdots,g_d)$ is actually in $\mathbb C^{\mu(Nx)+\tilde{\mu}(Nx)}.$, since by \ref{ortho1}, $h_i(x)=0$ if $i>\mu(Nx)$ and by \ref{ortho2}, $g_k(x)=0$ if $k>\tilde\mu(Nx)$. 

Before we generalize the Bratteli and Jorgensen definitions of $m$ system, we review the notions of Borel vector bundles over a Borel space and Borel cross-sections for Borel vector bundles.

\begin{definition} \cite{FD}
\label{bundledef}
Let $X$ be a topological space.  A (finite dimensional) vector bundle over the space $X,$ denoted by $(E,\;p,\;X),$ is a topological space $E,$ together with a continuous open surjection $p:\; E\;\rightarrow\;X,$ and operations and norms making each fiber $E_x\;=\;p^{-1}(X)$ into a (finite dimensional) vector space, which in addition satisfies the following conditions:
\renewcommand{\labelenumi}{(\roman{enumi})}
\begin{enumerate}
\item $y\mapsto\;\|y\|$ is continous from $E$ to $\mathbb R,$
\item The operation $+$ is continuous as a function from $\{(y,z)\in E\times E:\;p(y)=p(z)\}$ to $E.$
\item For each $\lambda\;\in\;\mathbb C,$ the map $y\mapsto\;\lambda\cdot y$ is continuous from $E$ to $E.$
\item If $x\in X$ and $\{y_i\}$ is any net of elements of $E$ such that $\|y_i\|\rightarrow 0$ and 
$p(y_i)\rightarrow x$ in $X,$ then $y_i\rightarrow \vec{0}\in E_x$ in $E.$ 
\end{enumerate}
A Borel map $s:X\;\rightarrow\; E$ is called a Borel cross-section if $p\circ s(x)\;=\;x,\forall x\in X.$
\end{definition}
Given an essentially bounded multiplicity function $\mu$ on $\mathbb T,$ let $c\;=\;\text{ess}\;\text{sup}\;\mu,\;d\;=\;\text{ess}\;\text{sup}\;\tilde{\mu},$ and let $T_j\;=\{x\in\mathbb T:\mu(Nx)+\tilde{\mu}(Nx)=j,\;0\;\leq\;j\;\leq\;c+d\}.$  Set $T_{i,j}\;=\;S_i\cap T_j,\;0\;\leq\;j\;\leq\;c+d;$ then each $T_{i,j}$ is Borel and $S_i\;=\;\bigsqcup_{j=0}^{c+d} T_{i,j}.$
\begin{definition}
\label{Msystemdef}
Let $E$ be the Borel space given by 
$$E\;=\;\bigsqcup_{i=1}^c \bigsqcup_{j=0}^{c+d} [T_{i,j}\times\mathbb C^j].$$
Let $(E,\;p,\;\bigsqcup_{i=1}^c S_i)$ be the Borel vector bundle where the map $p:\;E\;\rightarrow\;\:\bigsqcup_{i=1}^c S_i$ is defined by $p(x,\vec{v})=x,\;(x,\vec{v})\;\in\;T_{i,j}\times\mathbb C^j.$ By definition, $(E,\;p,\;\bigsqcup_{i=1}^c S_i)$ is a vector bundle $\bigsqcup_{i=1}^c S_i$ whose fiber over $x\in S_i$ is a complex vector space of dimension $\mu(Nx)+\tilde{\mu}(Nx).$ An {\bf $M$-system associated to the multiplicity function $\mu$} is a
Borel cross-section $M:\;\bigsqcup_{j=1}^c S_j\;\rightarrow\;E$ of this bundle whose values, $(M_1(x), M_2(x),\;\cdots M_{\mu(Nx)+\tilde{\mu}(Nx)})$, are the output of a vector of generalized low and high pass filters.
\end{definition}

Note all the information about the generalized filters $\{h_{i,j}\}$ and $\{g_{k,j}\}$ is encoded in the $M$-system  In particular, for any fixed multiplicity function $\mu$, such that both $\mu$ and $\tilde{\mu}$ are constant in a neighborhood of the origin, we have a one-to-one correspondence between M-systems and   
collections of generalized filter functions as defined in the previous section.  

To develop the loop group action on these $M$-systems, we start by defining an endomorphism $\Pi_N:\;\bigsqcup_{i=1}^c S_i\;\rightarrow\;\mathbb T$ by $\Pi_N(x)\;=\;Nx\;\text{mod}\;1.$ Each $x\in\mathbb T$ has $\sum_{l=0}^{N-1}\mu(\frac{x+l}{N})\;=\;\mu(x)\;+\;\tilde{\mu}(x)$ preimages in $\bigsqcup_{i=1}^c S_i.$
For convenience of notation, we label these preimage maps $r_{(l,j)},$ where $r_{(l,j)}(x)\;=\;\frac{x+l}{N}\;\in\;S_j\;\subseteq\;\bigsqcup_{i=1}^c S_i$ for 
$1\leq\;j\;\leq\mu(\frac{x+l}{N}).$   (Note that this range on $j$, as $l$ varies from 0 to $N-1$, gives all the preimages, since if $j>\mu(\frac{x+l}{N}),$ by definition $\frac{x+l}{N}$ is not an element of $S_j.$) For each fixed $x$, we give the pairs $(l,j)$ the lexicographical order, and thus implicitly define a 1-1 map $\lambda_x$ taking the pairs $(l,j)$ onto the integers from 1 to $\mu(x)+\tilde{\mu}(x)$. 

We now construct a unitary group bundle $(F,\; q,\;\mathbb T)$ as follows.  For each $j\;\in\;\{1,\cdots,c+d,\},$ let 
$Z_j\;=\;\{x\;\in\;\mathbb T:\;\mu(x)+\tilde{\mu(x)}=j\}.$ Let 
$$\widetilde{E}\;=\;\bigsqcup_{j=0}^{c+d} [Z_j\times M_j(\mathbb C)],$$
where $M_j(\mathbb C)$ is viewed as a $j^2$-dimensional normed vector space with the operator norm.
Defining $q:\;\widetilde{E}\;\rightarrow\;\mathbb T$ by $q(x, A)\;=\;=\;x),$ for $(x,A)\;\in\;E,$ we see that 
$(\widetilde{E},\;q,\mathbb T)$ is a Borel vector bundle in the sense of Definition \ref{bundledef}.
Now let $F$ be the subspace of $\widetilde{E}$ defined by 
$$F\;=\;\{(x,A)\in \widetilde{E}:\; A\in\;U(\mu(x)+\tilde{\mu}(x),\mathbb C)\}.$$
Then $q:\;F\;\rightarrow\;\mathbb T$ is a continuous open surjection, and the fiber $q^{-1}(x)$ of the bundle consists of the group of complex unitary matrices $U(\mu(x)+\tilde{\mu}(x),\mathbb C).$ Borel cross sections to this group bundle consist of Borel maps $K:\mathbb T\;\rightarrow\;F$ such that $q\circ K(x)\;=\;x.$
We denote the set of sections of this bundle by $\Gamma(F,q).$ Note $\Gamma(F,q)$ is a group under pointwise operations on $\mathbb T.$, where the identity element of the group is given by that section whose value at $x$ is equal to $Id_{\mu(x)+\tilde{\mu}(x)}.$ 

We are now ready to state a key theorem about $M$-systems that can be derived from the orthogonality relations:
\begin{theorem} 
\label{mainthm}
Let  $\Gamma(F,q)$ be the group of cross sections of the group bundle associated to the multiplicity function $\mu$ defined above.  Let $M:\;\bigsqcup_{j=1}^c S_j\;\rightarrow\;E$ be an $M$-system associated to $\mu.$ Then 
 $x\;\mapsto\;(K_{i,\lambda_x(l,j)}(x)),$ where 
$$K_{i,\lambda_x(l,j)}(x)\;=\;\sqrt{\frac{1}{N}}M_i(r_{(l,j)}(x))$$
is an element of $\Gamma(F,q).$
\end{theorem}
\begin{proof}
As noted above $1\leq\lambda_x(l,j)\leq\mu(x)+\tilde{\mu}(x)$, so for each $x\in\mathbb T,$ the matrix $(K_{i, \lambda_x(l,j)}(x))$ is a square matrix of the correct dimension.
We shall show that for all $x\in \mathbb T,$ the rows of $(K_{i,\lambda_x(l,j)}(x))$  are orthonormal. We use the orthogonality relations \ref{ortho1}, \ref{ortho2}, and \ref{ortho3} for generalized filter functions in this proof. 

Write $K_i$ for the $i$th row of $(K_{i,\lambda_x(l,j)}(x))$, and suppose first that $1\;\leq i\leq i'\leq\;\mu(x)$ Then 
\begin{eqnarray}
<K_i,\;K_{i'}>&=&\sum_{\lambda_x(l,j)=1}^{\mu(x)+\tilde{\mu}(x)}K_{i,\lambda_x(l,j)}(x)\overline{K_{i,\lambda_x(l,j)}(x)} \nonumber \\
&=& \sum_{l=0}^{N-1}\sum_{j=1}^{\mu(\frac{x+l}{N})}\sqrt{\frac{1}{N}}M_i(r_{(l,j)}(x))\overline{\sqrt{\frac{1}{N}}M_{i'}(r_{(l,j)}(x))} \nonumber \\
&=& \sum_{l=0}^{N-1}\sum_{j=1}^{\mu(\frac{x+l}{N})}\frac{1}{N}h_{i,j}(\frac{x+l}{N})\overline {h_{i',j}(\frac{x+l}{N})} \nonumber \\
&=& \frac{1}{N}\sum_{j=1}^{\mu(\frac{x+l}{N})}\sum_{l=0}^{N-1}h_{i,j}(\frac{x+l}{N})\overline {h_{i',j}(\frac{x+l}{N})} \nonumber \\
&=& \frac{1}{N}\sum_{j=1}^{c}\sum_{l=0}^{N-1}h_{i,j}(\frac{x+l}{N})\overline {h_{i',j}(\frac{x+l}{N})} \nonumber
\end{eqnarray}
(since $h_{i,j}(\frac{x+l}{N})\;=\;0$ for $j>\mu(\frac{x+l}{N})$ since $\frac{x+l}{N}\notin\;S_j$ in that case)  
$$\;=\;(\text{by}\;\ref{ortho1})\;\;N\frac{1}{N}\delta_{i,i'}\chi_{S_i}(x)\;=\;\delta_{i,i'}$$
(we note that $\chi_{S_i}(x)=1$ since we have $i\leq\;\mu(x)$ and for those values of $i,\;x\;\in\;S_i$ by definition of $\mu(x).$)

The cases $\mu(x)\;< i\leq i'\leq\;\mu(x)+\tilde{\mu}(x)$ and $1\;\leq i \leq\;\mu(x)<\; i'\leq\;\mu(x)+\tilde{\mu}(x)$  follow from similar arguments using \ref{ortho2} and \ref{ortho3}.  Thus we have that in all cases, the rows of $(K_{i,\lambda_x(l,j)}(x))$ are orthonormal, and we have the desired unitary matrix.
\end{proof}

The results of Theorem \ref{mainthm} imply that the columns of  $(K_{i,\lambda_x(l,j)}(x))$ are orthonormal as well, so we can deduce as a corollary: 
\begin{corollary} 
\label{orthcol}
Let $\mu$ and $\tilde{\mu}$ be  multiplicity and ``conjugate" multiplicity functions that are constant in a neighborhood of the origin, with related sequences of sets 
$\{S_i\:|\;1\leq\; i\;\leq c\}$ and $\{\tilde{S}_k\:|\;1\leq\; k\;\leq d\}.$ Suppose $\{h_{i,j}\}_{1\leq\;i,j\;\leq c}$ and 
$\{ g_{k,j}\}_{1\leq\;k\;\leq d,\;1\leq\;j\;\leq c}$ are generalized low-pass and high-pass filter functions with associated to the multiplicity function.  
Then for all $x\;\in\;\mathbb T,$ and for all $l,\;l'\;\in\;\{0,1,\cdots,N-1\}$ and $j,\;j'\;\in\;\{1,2,\cdots,c)\},$ we have 
$$\sum_{i=1}^{c+d}\frac{1}{N}M_i(r_{(l,j)}(x))\overline{M_i(r_{(l',j')}(x))}$$
$$\;=\;\sum_{i=1}^c \frac{1}{N}h_{i,j}(\frac{x+l}{N})\overline{h_{i,j'}(\frac{x+l'}{N})}\;+\;
\sum_{k=1}^d \frac{1}{N}g_{k,j}(\frac{x+l}{N})\overline{g_{k,j'}(\frac{x+l'}{N})}\;=\;
\delta_{j,j'}\delta_{l,l'},$$ 
where here the $M_i$ correspond to the generalized filter functions $h_{i,j}$ and $g_{k,j}$ as in Definition \ref{Msystemdef}. 
\end{corollary}
\begin{proof}
The statement that orthonormality of the rows implies orthonormality of the columns for finite-dimensional unitary matrices together with Theorem \ref{mainthm} tells us that 
for all $x\;\in\;\mathbb T,$ for all $l,l'\;\in\;\{0,1,\cdots,N-1\}$ and $j\;\in\;\{1,\cdots,\mu(\frac{x+l}{N})\},\;j'\;\in\;\{1,\cdots,\mu(\frac{x+l'}{N})\}$ we have 
$$\sum_{i=1}^{\mu(x)+\tilde{\mu}(x)}\frac{1}{N}M_i(r_{(l,j)}(x))\overline{M_i(r_{(l',j')}(x))}$$
$$\;=\;\sum_{i=1}^{\mu(x)}\frac{1}{N}h_{i,j}(\frac{x+l}{N})\overline{h_{i,j'}(\frac{x+l'}{N})}\;+\;
\sum_{k=1}^{\tilde{\mu}(x)}\frac{1}{N}g_{k,j}(\frac{x+l}{N})\overline{g_{k,j'}(\frac{x+l'}{N})}\;=\;
\delta_{j,j'}\delta_{l,l'}.$$
But now we recall Equation \ref{ortho1}, which tells us that for fixed $x\in\;\mathbb T,$ if $i$ is a fixed integer satisfying $\mu(x)\;<\; i\;\leq\; c,$ then 
$$\sum_{j=1}^c\sum_{l=0}^{N-1}|h_{i,j}(\frac{x+l}{N})|^2\;=\;0$$ (since $x\;\notin\; S_i$ for $i>\mu(x)$ by definition of $\mu(x).)$  Thus we must have
$$h_{i,j}(\frac{x+l}{N})\;=\;0,\mu(x)\;<\; i\;\leq\; c,\;1\;\leq\; j\;\leq c,\;0\;\leq\; l\;\leq N-1.$$
Similarly, equation \ref{ortho2} tells us that for fixed $x\;\in\;\mathbb T,$ if $k$ satisfies $\tilde{\mu(x)}\;<\; k\;\leq\; d,$ then 
$$\sum_{j=1}^c\sum_{l=0}^{N-1}|g_{k,j}(\frac{x+l}{N})|^2\;=\;0$$ (since $x\;\notin\;\widetilde{S_k}$ for $k>\tilde{\mu}(x),$ by definition of the $\widetilde{S_k}).$  Thus we must have
$$g_{k,j}(\frac{x+l}{N})\;=\;0,\;\tilde{\mu}(x)\;<\; k\;\leq\; d,\;1\;\leq\; j\;\leq c,\;0\;\leq\; l\;\leq N-1.$$
Finally we note that whenever $j> \mu(\frac{x+l}{N}),$ we will have the identities $h_{i,j}(\frac{x+l}{N})\;=\;g_{k,j}(\frac{x+l}{N})\;=\;0$ for all $i$ between $1$ and $c$ and for all $k$ between $1$ and $d.$ The reason for this is that in this case, $\frac{x+l}{N}\notin\; S_j$ and for fixed $j,$ all the functions $\{h_{i,j},\;g_{k,j}\}$ are supported on $S_j.$ From this we see that for all $x\;\in\;\mathbb T,$ and for all $l,\;l'\;\in\;\{0,1,\cdots,N-1\}$ and $j,\;j'\;\in\;\{1,\cdots,\mu(\frac{x+l}{N})\},$ we have 
$$\sum_{i=1}^c \frac{1}{N}h_{i,j}(\frac{x+l}{N})\overline{h_{i,j'}(\frac{x+l'}{N})}\;+\;
\sum_{k=1}^d \frac{1}{N}g_{k,j}(\frac{x+l}{N})\overline{g_{k,j'}(\frac{x+l'}{N})}\;=\;
\delta_{j,j'}\delta_{l,l'},$$
since the extra terms we are adding on to the left hand side not coming from the orthonormality of columns are all zero, by the remarks above.
\end{proof}

The following example explicitly shows how Theorem \ref{mainthm} works in the particular case of the Journ\'e wavelet.  The construction of the generalized filters for the Journ\'e wavelet was first done in the thesis of J. Courter \cite{CO}. 
\begin{example}
The Journ\'e wavelet in the frequency domain is 
the characteristic function of the set 
$$[-\frac{16}{7},-2)\cup [-\frac{1}{2},-\frac{2}{7})\cup [\frac{2}{7},\frac
{1}{2}]\cup [2,\frac{16}{7}).$$
Here $\mu$ takes on the values $0,1,$ and $2,$ and $\widetilde{\mu(x)}\equiv 1,$ since the Journ\'e wavelet is a single orthonormal wavelet.
If we identify $\mathbb T$ with  $[-\frac{1}{2},\frac{1}{2}),$ we can write  
$S_1=[-\frac {1}{2},-\frac{3}{7})\cup [-\frac {2}{7},\frac {2}{7})\cup 
[\frac {3}{7},\frac{1}{2}),\;S_2=[-\frac {1}{7},\frac{1}{7}),$ and 
$\widetilde{S_1}=[-\frac {1}{2},\frac{1}{2}].$
The generalized filter functions then are: 
$$h_{1,1}(x)=\chi_{[-\frac{2}{7},-\frac{1}{4})\cup (-\frac{1}{7},\frac{1}{7})\cup 
[\frac{1}{4},\frac{2}{7})}(x),$$
$$h_{1,2}(x)=0,$$
$$h_{2,1}=\chi_{[-\frac{1}{2},-\frac{3}{7})\cup [\frac{3}{7},\frac{1}{2})}(x),$$
$$h_{2,2}(x)=0;$$
$$g_1(x)=\chi_{[-\frac{1}{4},-\frac{1}{7})\cup [\frac{1}{7},\frac{1}{4})}(x),$$
$$g_2(x)=\chi_{[-\frac{1}{7},\frac{1}{7})}(x).$$
\noindent Consider the decomposition of the circle $\mathbb T$ (identified with $[-\frac{1}{2},\frac{1}{2})$) 
given by
\vskip.1in

$P_1=[-\frac{1}{7},\frac{1}{7})$ (Here $\mu(x)=2,\mu(\frac{x}{2})=2,\mu(\frac{x+1}{2})=1.)$

$P_2=\pm [\frac{1}{7},\frac{2}{7})$ (Here $\mu(x)=1,$ $\mu(\frac{x}{2})=2,$ $\mu(\frac{x+1}{2})=0.)$

$P_3=\pm [\frac{2}{7},\frac{3}{7})$, (Here $\mu(x)=0$, $\mu(\frac{x}{2})=1,\;\mu(\frac {x+1}2)=0.$)

$P_4=\pm [\frac{3}{7},\frac{1}{2})$  (Here $\mu(x)=1,$ $\mu(\frac{x}{2})=1,\;\mu(\frac {x+1}2)=1.$)
\vskip.1in
\noindent The associated cross-section matrix bundle is:
$$K_{i,\lambda_x(l,j)}(x)=\left\{\begin{array}{ll}
{\left(\begin{array}{ccc}
1\;&0\;&0\\
0\;&0\;&1\\
0\;&1\;&0\end{array}\right),}&\mbox{if}\ \;x\in P_1,\\
{\left(\begin{array}{cc}
1\;&0\\
0\;&1\end{array}\right),}&\mbox{if}\ \;x\in P_2,\\
\;\;\;\;\;\;\;\;1,&\mbox{if}\ \;x\in P_3,\\
{\left(\begin{array}{cc}
0\;&1\\
1\;&0\end{array}\right),}&\mbox{if}\;x\in P_4.
\end{array}\right.$$
\end{example}

We are ready to define the generalized loop group and its associated action on the set of $M$-systems associated to a multiplicity function $\mu$.
\begin{definition} 
The {\bf loop group} associated to the multiplicity function $\mu$ 
is defined to be the subgroup $\text{Loop}(F,q)$ of the group of Borel sections $\Gamma(F,q)$ whose elements $K$ 
satisfy $K(0)\;=\;Id_{\mu(0)+\tilde{\mu}(0)},$ and  $K_{i,\lambda_x(j,l)}$ are $\text{Lip}_1$ in a neighborhood of the origin.  
\end{definition}

We now come to the main theorem of the paper.  Just as in the classical case, it is possible to show that the generalized loop group acts freely and transitively on the set of $M$-systems:
\begin{theorem}
There is a free and transitive action of $\text{Loop}(F,q)$ on the set of $M$-systems associated to an essentially bounded multiplicity function $\mu$ such that $\mu$ is constant in neighborhoods of $\frac{l}{N},\;0\leq l\leq N-1.$ This action is given by 
$$K\cdot M(x)\; =\; K(\Pi_N (x))[(M_1(x),M_2(x),\cdots, M_{\mu(Nx)+\tilde{\mu}(Nx)}(x))]^T.$$
\end{theorem}
\begin{proof}
We prove the transitivity first.  Suppose we are given two different $M$-systems, labeled $M=(M_i)$ and $\widetilde{M}=(\widetilde{M}_i).$ Define an element $K$ of the group bundle, that is, an element of $\Gamma(F,q),$ 
where $K(x)$ has dimension $\mu(x)+\tilde {\mu}(x)$, as follows:
$$K_{i,i'}(x)=\frac{1}{N}\sum_{\lambda_x(l,j)=1}^{\mu(x)+\tilde{\mu}(x)} \overline {M_i'(r_{(l,j)}(x))}\widetilde{M_i}(r_{(l,j)}(x)).$$  
As in the classical case, we have 
\begin {eqnarray}
[K\cdot M]_i(x)&=& \sum_{i'=1}^{\mu(Nx)+\tilde{\mu}(Nx)}K_{i,i'}(\Pi_N (x))M_i'(x) \nonumber \\
&=& \sum_{i'=1}^{\mu(Nx)+\tilde{\mu}(Nx)}\left(\frac{1}{N}\sum_{\lambda_{Nx}(l,j)=1}^{\mu(Nx)+\tilde{\mu}(Nx)} \overline {M_i'(r_{(l,j)}(Nx))}\widetilde{M_i}(r_{(l,j)}(Nx))\right)M_i'(x) \nonumber \\
&=& \sum_{\lambda_{Nx}(l,j)=1}^{\mu(Nx)+\tilde{\mu}(Nx)}\widetilde{M_i}(r_{(l,j)}(Nx))\left(\sum_{i'=1}^{\mu(Nx)+\tilde{\mu}(Nx)}\overline {M_i'(r_{(l,j)}(Nx))}M_i'(x)\right) \nonumber \\
&=& \widetilde{M_i}(x), \nonumber
\end {eqnarray}
where the last equality follows since, by the orthogonality of the columns of $M$ as established in Corollary \ref{orthcol}, we have that inside sum is 0 except for the single values of $l$ and $j$ where $r_{l,j}(Nx)=x$.  

Now to establish that the action is free, suppose $M=(M_i)$ is a $M$-system associated to $\mu$ and $K\;\in\;\text{Loop}(F,q)$ satisfies 
$$K(\Pi_N(x))[(M_1(x),M_2(x),\cdots, M_{\mu(Nx)+\tilde{\mu}(Nx)}(x))]^T\;=$$
$$[(M_1(x),M_2(x),\cdots, M_{\mu(Nx)+\tilde{\mu}(Nx)}(x))]^T.$$
Then, analogously to the classical case, for each $x\in\mathbb T$ we define a $(\mu(Nx)+\tilde{\mu}(Nx))\times (\mu(Nx)+\tilde{\mu}(Nx))$ unitary matrix ${\mathcal M}$ by 
$${\mathcal M}_{i,\lambda_{Nx}(l,j)}(x)\;=\;\sqrt{\frac{1}{N}}M_i (r_{(l,j)}(Nx)),$$
$$1\;\leq\; i\;\leq\; \mu(Nx)+\tilde{\mu}(Nx),\;0\leq\;l\;\leq\;N-1,\;1\leq\;j\;\leq\mu(x+\frac{l}{N}).$$  
We then see that $K(Nx){\mathcal M}(x)\;=\;{\mathcal M}(x)$ for all $x\in\;\mathbb T.$ 
By unitarity of ${\mathcal M}(x),$ this shows that  that $K(Nx)$ is the $(\mu(Nx)+\tilde{\mu}(Nx))\times (\mu(Nx)+\tilde{\mu}(Nx))$ identity matrix for all $x\in\;\mathbb T,$ which implies that $K$ is the identity element of $\text{Loop}(F,q),$ as desired. 
\end{proof}

We note this result is set purely in the language of transformation groups, i.e. we have defined a set corresponding to a fixed multiplicity function (the $M$-systems), and we have described a group which acts freely and transitively on this set.  In a later paper we intend to impose conditions on $M$-systems under which we can mimic the classical construction and use them to obtain a normalized tight frame wavelet family for dilation by $N.$

\begin{acknowledgments}
The authors gratefully acknowledge helpful suggestions from 
Astrid An Huef and Iain Raeburn.
They also thank Brian Treadway for coordinating and merging \TeX\ files.
\end{acknowledgments}

\bibliographystyle{amsalpha}

\end{document}